\documentclass{article}
\usepackage{amsmath,amsfonts,amsthm}

\providecommand{\C}[1]{\mathcal{#1}}
\providecommand{\D}[1]{\mathbb{#1}}

\newcommand{\ee}{\mathrm{e}}

\def\K1{{{\cal {A}}_1}}
\def\k{\kappa}

\newcommand\capp{{\rm cap}}

\newcommand\sy[1]{(\ref{#1})}
\theoremstyle{plain}
\newtheorem{thm}{Theorem}%[section]

\newtheorem{prop}[thm]{Proposition}

\newtheorem{rem}[thm]{Remark}
\begin{document}
\title{Capacity of a multiply-connected domain and nonexistence of
  Ginzburg-Landau minimizers with prescribed degrees on the boundary.}
%
%\author{Leonid Berlyand\inst{1} \and Dmitry Golovaty\inst{2} \and
%  Volodymyr Rybalko\inst{3}\thanks{The first author was supported by
%    the NSF grant DMS-0204637, the second author was supported by the
%    NSF grant DMS-0407361, and the third author was supported by the
%    grant GP/F8/0045 $\sf \Phi$8/308-2004}}
%
\author{L.~Berlyand \thanks{Department of Mathematics, The
    Pennsylvania State University, University Park, PA 16802, USA,
    email: {\em{berlyand@math.psu.edu}}. Supported in part by the NSF
    grant DMS-0204637.} \and D.~Golovaty \thanks{Department of
    Theoretical and Applied Mathematics, The University of Akron,
    Akron, OH 44325, USA, email: {\em{dmitry@math.uakron.edu}}.
    Supported in part by the NSF grant DMS-0407361.} \and V.~Rybalko
  \thanks{Mathematical Division, Institute for Lower Temperature
    Physics and Engineering, 47 Lenin Ave., 61164 Kharkov, Ukraine,
    email: {\em{vrybalko@ilt.kharkov.ua}}. Supported in part by the
    grant GP/F8/0045 $\sf \Phi$8/308-2004.}}
\maketitle
\date
\begin{abstract}
  Let $\omega$ and $\Omega$ be bounded simply connected domains in
  ${\D R}^2$, and let $\bar\omega\subset\Omega$. In the annular domain
  $A=\Omega\setminus\bar\omega$ we consider the class $\C J$ of
  complex valued maps having modulus $1$ and degree $1$ on $\partial
  \Omega$ and $\partial\omega$.
  
  It was conjectured in \cite{bemi04a} that the existence of
  minimizers of the Ginzburg-Landau energy $E_\k$ in $\C J$ is
  completely determined by the value of the $H^1$-capacity $\capp(A)$
  of the domain and the value of the Ginzburg-Landau parameter $\k$.
  The existence of minimizers of $E_\k$ for all $\k$ when
  $\capp(A)\geq\pi$ (domain $A$ is ``thin'') and for small $\k$ when
  $\capp(A)<\pi$ (domain $A$ is ``thick'') was demonstrated in
  \cite{bemi04a}.
  
  Here we provide the answer for the case that was left open in
  \cite{bemi04a}. We prove that, when $\capp(A)<\pi$, there exists a
  {\em{finite}} threshold value $\k_1$ of the Ginzburg-Landau
  parameter $\k$ such that the minimum of the Ginzburg-Landau energy
  $E_\k$ {\em{not}} attained in $\C J$ when $\k>\k_1$ while it is
  attained when $\k<\k_1$.
\end{abstract}
\section{Introduction}
The present paper establishes nonexistence of minimizers of the
Ginzburg-Landau functional in a class of Sobolev functions with
prescribed degree on the boundary of an annular domain when the
$H^1$-capacity of the domain is less than the critical value
$c_{cr}=\pi$. Here an annular domain is any domain in ${\D R}^2$
conformal to a circular annulus.

\subsection{Mathematical formulation and physical model}
Consider the minimization problem for the Ginzburg-Landau functional
\begin{equation}
E_\kappa[u]=\frac{1}{2}\int_A|\nabla u|^2dx
+\frac{\kappa^2}{4}\int_A (|u|^2-1)^2dx\to {\rm inf},
\quad
u\in {\C J},
\label{P1}
\end{equation}
where $A=\Omega\setminus\bar\omega$, $\bar\omega\subset \Omega$, and
$\omega,\ \Omega$ are bounded, simply connected domains in ${\D R}^2$
with smooth boundaries. The class ${\C J}$ is defined by
\begin{equation}
{\C J}=\{u\in H^1(A):
|u|=1 \ \text{on}\ \partial \Omega\cup\partial\omega,\ 
{\rm deg}(u,\partial \Omega)={\rm deg}(u,\partial\omega)=1\}.
\label{classJ}
\end{equation}
Note that a minimizer of \sy{P1} in $\C J$ satisfies the
Ginzburg-Landau equation
\begin{equation}
-\Delta u+\kappa^2(|u|^2-1)u=0\,,
\label{(0)}
\end{equation}
in $A$ along with the natural boundary conditions \( \frac{\partial
  u}{\partial \nu}\times u=0\quad \text{on}\ \partial A.  \)

The problem \sy{P1} originates with a Ginzburg-Landau variational
model of superconducting persistent currents in multiply connected
domains.

Consider a superconducting material with a hole occupying the domain
$A$. The superconductor is characterized by a complex order parameter
$u$ with the magnitude $|u|$ and the gradient of $\mbox{arg}\, u$
describing the density of superconducting electrons and
superconducting current, respectively. The order parameter vanishes in
a normal, non-superconducting state while it is $S^1$-valued in a
perfectly superconducting state.

In a ``full'' superconductivity problem, the Ginzburg-Landau
functional depends not only on the order parameter but also on
magnetic field. In order to facilitate theoretical analysis, various
simplifications of the Ginzburg-Landau functional have been introduced
and studied \cite{Abrikosov}, \cite{bbh94}.

From now on, suppose that the external magnetic field is zero and a
characteristic size of the domain $A$ is smaller than the penetration
depth. Then the current-induced magnetic field can be neglected and
the energy of the superconductor reduces to the functional
in \sy{P1}. In the absence of the external field, the {\em{local}}
minimizers of $E_\kappa$ in $H^1$ can be interpreted as
{\em{persistent currents}} in a superconducting composite with holes
\cite{deaver}.

Alternatively, the same currents can be understood as {\em{global
    minimizers}} of the Ginzburg-Landau functional $E_\kappa$ over $A$
when the order parameter $u$ is in the class of complex-valued maps
with a prescribed degree on each connected component of the boundary.
The degree boundary condition reflects the topological quantization of
the phase of the order parameter around the hole.

\subsection{Existing results}
The Ginzburg-Landau-theory-related literature is too vast a topic to
be fully explored within the limited scope of this paper and we will
restrict our review to studies that are most relevant to our problem.

The asymptotics as $\kappa\to \infty$ of global minimizers for the
Ginzburg-Landau functional and their vortex structure for the
Dirichlet boundary data (for which the degree is fixed by default)
were studied in detail in \cite{bbh94} for simply-connected domains.

A minimization problem for the Ginzburg-Landau functional with a
magnetic field for classes of functions with no prescribed boundary
conditions in simply connected domains was studied in
\cite{serfaty1}-\nocite{serfaty4}\cite{serfaty3}.  In this case, the
qualitative changes in the behavior of minimizers are governed by the
magnitude of the external magnetic field. In particular, the existence
of a threshold field value corresponding to a transition from
vortex-less minimizers to minimizers with vortices was proved in
\cite{serfaty4} when $\kappa\to\infty$.
 
The existence of local minimizers of the Ginzburg-Landau functional
with a magnetic field over three-dimensional tori was considered in
\cite{rubsternberg} (see \cite{jimbomorita2} for related results for
solids of revolution with a convex cross-section). The approach of
\cite{rubsternberg} relies on the fact that, when the parameter
$\kappa$ is large, the boundedness of the nonlinear term in the
Ginzburg-Landau energy forces the minimizing maps to be ``close'' to
$S^1$-valued maps.  The first step in the proof consisted of finding
local minimizers for the Dirichlet integral in all homotopy classes of
$S^1$-valued maps. Then, for $\kappa$ large, the existence of a local
minimizer of the Ginzburg-Landau functional in a vicinity of each
minimizer of the Dirichlet integral was shown.

The existence and properties of {\em{global}} minimizers of the
Ginzburg-Landau functional describing a superconductor in the presence
of magnetic field was studied in \cite{AB04} over multiply connected
domains. In the singular limit of $\kappa\to\infty$, it was
established that the holes in the domain act as ``giant vortices''
when the external field is fixed independent of $\kappa$. Further,
when the external field is of order $\log{\kappa}$, the interior
vortices start to appear in the domain, resembling the results of
\cite{serfaty1}-\nocite{serfaty4}\cite{serfaty3}.

In a related problem, the global minimizers of the Ginzburg-Landau
functional describing the uniformly rotating Bose-Einstein condensate
in a circular domain were considered in \cite{afetal04}.  Although the
domain in \cite{afetal04} is simply connected, in the limit the
solution is effectively restricted to the doubly-connected domain (an
annulus) by assuming that the pinning term vanishes in a smaller
circular region centered at the origin.

Note that, for none of the results mentioned so far, the existence and
the qualitative behavior of minimizers depend on the $H^1-$capacity
\cite{bemi04} of the domain
\begin{equation}
\label{capacity}
\mbox{\rm cap}(A)=\mbox{\rm Min }\left\{\int_A|\nabla v|^2\ ;\ v\in H^1(A), v=0\ \mbox{\rm on }\partial\Omega, v=1\ \mbox{\rm on }\partial\omega\right\}.
\end{equation}
 
The questions of existence and uniqueness of minimizers of the
Ginzburg-Landau functional in a class of maps with the degree boundary
conditions were studied in \cite{bervoss}, \cite{BG2002},
\cite{bemi04a}-\nocite{bemi04}\cite{bemi03}. For a narrow circular
annulus both existence and uniqueness were proved in \cite{BG2002} for
an {\em{arbitrary}} (not necessarily large) $\kappa>0$.  The
techniques of \cite{BG2002} rely on a priori estimates valid for
radially symmetric domains and cannot be readily extended to arbitrary
multiply-connected domains. In
\cite{bemi04a}-\nocite{bemi04}\cite{bemi03} a general approach for
such domains was developed. It was shown in
\cite{bemi04a}-\nocite{bemi04}\cite{bemi03} that, when the capacity of
a domain exceeds a certain critical value, the global minimizers of
the Ginzburg-Landau functional with degree boundary conditions exist
for arbitrary $\kappa$. These minimizers are vortexless and unique for
large $\kappa$.  When the capacity is below the critical value, the
minimizing sequences must develop vortices near the boundary of the
domain for large $\kappa$. When the domain is conformally equivalent
to an annulus and the degrees of admissible functions are equal to $1$
on both connected components of the boundary, it was proved in
\cite{bemi04} that minimizing sequences develop {\em{exactly}} two
vortices of degree $1$ and $-1$.

Mathematically, the assumption that $|u|=1$ on the boundary has a very
interesting implication that the vortices in domains of small capacity
can approach the boundary at least exponentially close in $\kappa$.
Consequently, it has been exceedingly difficult to demonstrate whether
the vortices in a minimizing sequence actually end up on the boundary
itself or they reach a limiting point in the interior of the domain.

The existence of the critical domain capacity for functions with
degree boundary conditions is related to the fact that the class
${\mathcal J}$ is not closed with respect to weak $H^1$-topology
\cite{bemi03}. Further, the attainability of the lower bound for
$E_\kappa[u]$ over ${\mathcal J}$ cannot be deduced by the direct
method of calculus of variations.

Recall the following results from \cite{bemi03}.

\begin{thm}
Assume that $\capp(A)\geq\pi$. Then 
\begin{equation}
m_\kappa={\rm Inf}\left\{E_\kappa[u],u\in{\C J}\right\}
\label{minimum}
\end{equation}  
is attained for all $\kappa>0$.
\end{thm}

\begin{thm}
\label{thm2}
Assume that $\capp(A)<\pi$. Then either $m_\kappa$ is attained for all
$\kappa>0$ or there exists a $\kappa_1<\infty$ such that $m_\kappa$ is
always attained for $\kappa<\kappa_1$ and it is never attained for
$\kappa>\kappa_1$.
\end{thm}

It was conjectured in \cite{bemi04a}-\nocite{bemi04}\cite{bemi03} that
the second of the two cases in Theorem \ref{thm2} always occurs, that
is there is a threshold value of $\kappa_1$ above which {\em{the
    minimizer does not exist}} in supercritical domains. The existence
of $\kappa_1$ is established in this paper.

\subsection{Main result and outline of the proof}
Our main result is the following

\begin{thm} 
\label{mainth}
Assume $\capp(A)<\pi$. Then there is a finite $\kappa_1>0$ such that
$m_\kappa$ is always attained for $\kappa<\kappa_1$ and it is never
attained for $\kappa>\kappa_1$.
\end{thm}

The proof of Theorem \ref{mainth} is based on the estimate
$m_\kappa\leq 2\pi$ established in \cite{bervoss} as well as on the
convergence results in \cite{bemi03}.

We argue by contradiction. Assume that there is no finite $\kappa_1$
such that $m_{\kappa_1}$ is not attained (or, equivalently,
$\kappa_1=\infty$). Then the minimizer of $E_\kappa$ exists for all
finite $\kappa$.

Now suppose that ${\C A}$ is a circular annulus conformally equivalent
to $A$. Then $\capp({\C A})=\capp(A)$. As we show in Section
\ref{conformaleq}, given our assumptions in the previous paragraph, the
minimum of $E_\kappa$ over ${\C A}$ is attained for all finite
$\kappa$ and one can assume without loss of generality that $A$ is a
circular annulus.

Since $m_\kappa$ is attained for all $\kappa$, for every $\kappa>0$
there exists a $u_\kappa\in{\C J}$ such that
$E\left[u_\kappa\right]=m_\kappa\leq 2\pi$.

Next, in Section \ref{construction} we construct a sequence of
auxiliary quadratic functionals $\left\{F_\kappa\right\}_{\kappa>0}$
over a rectangular domain with linear Euler-Langrange equations and
use $\left\{u_\kappa\right\}_{\kappa>0}$ to produce a sequence of
functions $\left\{v_\kappa\right\}_{\kappa>0}$ such that
$F_\kappa\left[v_\kappa\right]\leq 2\pi$.

Finally, we complete the proof in Section \ref{estimate} by finding
the explicit solution $w_\kappa$ of the system of linear PDEs
corresponding to $F_\kappa$ and use this solution to show that
$F_\kappa\left[w_\kappa\right]>2\pi$.

\section{Preliminary results}
\label{prelims}
Here we gather prior results from \cite{bemi03}, \cite{bemi04a}, and
\cite{bervoss} that will be needed to prove Theorem \ref{mainth}.

\begin{prop}{\rm{(\cite{bemi03})}}
  Assume that $m_\kappa<2\pi$. Then $m_\kappa$ is attained.
\label{exists}
\end{prop}

The bound $2\pi$ for $m_\kappa$ is, in fact, precise due to
to

\begin{prop}{\rm(\cite{bervoss})}
\label{bound2pi}
For all $\kappa>0$ we have
$$
m_\kappa\leq2\pi.
$$ 
\end{prop}

Finally, recall the following theorem from \cite{bemi04a}.

\begin{thm} {\rm(\cite{bemi04a})}
  Let $\capp(A)<\pi$, and suppose that $u_\kappa\in {\C J}$ is a
  solution of Ginzburg-Landau equation {\em{(\ref{(0)})}} such that
  $E_\kappa(u_\kappa)<2\pi+\ee^{-\kappa}$. Then there is
  $\gamma_\kappa={\rm const}\in S^1$ such that for any compact set $K$
  in $A$
\begin{equation}
\|u_\kappa-\gamma_\kappa\|_{C^l(K)}=o(\kappa^{-m}),\
\text{as}\ \kappa\to\infty,\ \forall m>0,\ l\in\D{N},
\label{(10)}
\end{equation}
\begin{equation}
\int_A(|u_\kappa|^2-1)^2 dx=o(\kappa^{-m}),\
\text{as}\ \kappa\to\infty,\ \forall m>0.
\label{(11)}
\end{equation}
\end{thm}

\section{Proof of Theorem \ref{mainth}}
We argue by contradiction. Suppose that for all $\kappa>0$, the infimum
$m_\kappa$ is attained at some map $u_\kappa\in\C J$.  Then, in view
of Proposition \ref{bound2pi},
\[E_\kappa\left[u_\kappa\right]\leq2\pi.\]
Next, we show that, without loss of generality, we can assume that $A$
is a circular annulus $A=\{x\in\D{R}^2: R>|x|>\frac{1}{R}\}$.

\subsection{Conformal equivalence to a circular annulus}
\label{conformaleq}
\begin{prop}
\label{conformal}
Suppose that $A$ is such that $m_\kappa$ is attained for every
$\kappa>0$. Then the same holds for the annular domain $${\C A}:=
\left\{x:\exp{\left(-\frac{\pi}{{\rm cap}(A)}\right)}<|x|<
  \exp{\left(\frac{\pi}{{\rm cap}(A)}\right)}\right\}\,,$$
where ${\C
  A}$ is conformally equivalent to $A$.
\end{prop}

\begin{proof} 
  First, observe (\cite{complanalysis}) that $A$ is conformally
  equivalent to a circular annulus ${\C A}$; moreover the
  corresponding conformal map ${\C F}$ extends to a
  $C^1$-diffeomorphism of $\bar A$ onto $\bar {\C A}$ that preserves
  the orientation of curves.
  
  Let $u_\kappa$ ($\kappa>0$) be a minimizer of the functional
  $E_\kappa[u]$ in ${\C J}$, then
\begin{equation}
  m_\kappa= E_\kappa\left[u_\kappa\right]<2\pi.
\label{bound} 
\end{equation}
Indeed, for any $\kappa^\prime>\kappa$ there is a minimizer
$u_{\kappa^\prime}$ of $E_{\kappa^\prime}[u]$ in ${\C J}$ and
$$
E_{\kappa^\prime}\left[u_{\kappa^\prime}\right]
\leq 2\pi,
$$
by Proposition \ref{bound2pi}. Then
$$
E_{\kappa^\prime}\left[u_{\kappa^\prime}\right]-E_\kappa\left[u_\kappa\right]\geq
E_{\kappa^\prime}\left[u_{\kappa^\prime}\right]-E_\kappa\left[u_{\kappa^\prime}\right]=
\frac{(\kappa^\prime)^2-\kappa^2}{4}
\int_{A}(|u_{\kappa^\prime}|^2-1)^2dx,
$$
so that $m_\kappa\leq 2\pi$ and $m_\kappa=2\pi$ if and only if
$|u_{\kappa^\prime}|=1$ a.e. in $A$. The map $u_{\kappa^\prime}$ is a
solution of Ginzburg-Landau equation (\ref{(0)}) because
$u_{\kappa^\prime}$ minimizes $E_{\kappa^\prime}[u]$ with respect to
its own boundary data. The pointwise equality $|u_{\kappa^\prime}|=1$
a.e. in $A$ implies that the phase of $u_{\kappa^\prime}$ satisfies
the Laplace equation in $A$ subject to the homogeneous Neumann
boundary conditions on $\partial A$. Then $u_{\kappa^\prime}\equiv{\rm
  const}$, in contradiction with $u_{\kappa^\prime}\in{\C J}$ and we
arrive at (\ref{bound}).

By using the conformal change of variables $x\to{\C F}(x)$, we obtain
from (\ref{bound}) that
$$
\frac{1}{2}\int_{\C A}|\nabla \tilde u|^2dx+ \frac{\kappa^2}{4}
\int_{\C A}(|\tilde u|^2-1)^2{\rm Jac}({\C F}^{-1})dx<2\pi,
$$
where $\tilde u(x)=u_\kappa({\C F}^{-1}(x))$. Since $\kappa$ is
arbitrary, using Proposition \ref{exists} we obtain that the minimum
of (\ref{P1}) is attained for all $\kappa>0$.
\end{proof}

\begin{rem}
  \label{remark1}
  Suppose that ${\C A}$ is as defined in Proposition \ref{conformal}.
  By \sy{capacity} and conformal invariance of the Dirichlet integral,
  we have that $\capp\,({\C A})=\capp\,(A)$.
\end{rem}

\subsection{Proof of Theorem \ref{mainth} continued: reduction to a linear problem}
\label{construction}
Multiplying the equation (\ref{(0)}) by $\log \frac{|x|}{R}$ and
integrating over $D=\{x:1<|x|<R\}$ we obtain
$$
0=\int_{D}\Delta u_\kappa \log \frac{|x|}{R}\,dx+
\kappa^2\int_{D}u_\kappa(1-|u_\kappa|^2)\log \frac{|x|}{R}dx
$$
$$
= \int_{\partial D} \frac{\partial u_\kappa}{\partial \nu} \log
\frac{|x|}{R} d\sigma -\int_{\partial D} u_\kappa\frac{\partial
  \log|x|}{\partial \nu} d\sigma
+\kappa^2\int_{D}u_\kappa(1-|u_\kappa|^2)\log \frac{|x|}{R}dx
$$
$$
=-\frac{1}{R}\int_{|x|=R}u_\kappa d\sigma+ \int_{|x|=1}u_\kappa
d\sigma
$$
$$+\int_{|x|=1} \frac{\partial u_\kappa}{\partial \nu} \log
\frac{1}{R} d\sigma+ \kappa^2\int_{D}u_\kappa(1-|u_\kappa|^2)\log
\frac{|x|}{R}dx.
$$
Therefore, by using (\ref{(10)}) and (\ref{(11)}) we have, as
$\kappa\to\infty$,
\begin{equation}
\frac{1}{R}\int_{|x|=R}u_\kappa d\sigma=2\pi\gamma_\kappa+o(\kappa^{-m}).
\label{(12)}
\end{equation}
A similar calculation over $D=\{x:R^{-1}<|x|<1\}$ leads to the
estimate
\begin{equation}
R\int_{|x|=1/R}u_\kappa d\sigma=2\pi\gamma_\kappa+o(\kappa^{-m}).
\label{(13)}
\end{equation}

Changing the variables $x\to (r,\varphi): \ x=e^{r+i\varphi}$, we have
$$
E_\kappa\left[u_\kappa\right]=\frac{1}{2}\int_{-L}^L
dr\int_{0}^{2\pi}d\varphi |\nabla u_\kappa|^2 +
\frac{\kappa^2}{4}\int_{-L}^L e^{2r} dr\int_{0}^{2\pi}d\varphi
(|u_\kappa|^2-1)^2,
$$
where $-{\rm log} R<r<{\rm log} R$, $0\leq\varphi<2\pi$, and
$L={\rm log} R$. 

We modify $u_\kappa$ as follows. First, to simplify the subsequent
calculations, we set either
$$
u^{(1)}_\kappa(r,\varphi):=\bar\gamma_\kappa
\begin{cases}
  u_\kappa(r,\varphi), & 0\leq r<L,\\
  u_\kappa(-r,\varphi), & -L<r<0,
\end{cases}
$$
or
$$
u^{(1)}_\kappa(r,\varphi):=\bar\gamma_\kappa
\begin{cases}
  u_\kappa(-r,\varphi), & 0\leq r<L,\\
  u_\kappa(r,\varphi), & -L<r<0,
\end{cases}
$$
to obtain that $u^{(1)}_\kappa(r,\varphi)=u^{(1)}_\kappa(-r,\varphi)$ and
\begin{equation}
\frac{1}{2}\int_{-L}^L dr\int_{0}^{2\pi}d\varphi
|\nabla u^{(1)}_\kappa|^2 +
\frac{\kappa^2}{4R^2}\int_{-L}^L  dr\int_{0}^{2\pi}d\varphi
\left(|u^{(1)}_\kappa|^2-1\right)^2\leq 2\pi\,.
\label{(14)}
\end{equation}
Due to (\ref{(10)}), for all $0<\rho <L$ and $m\in\D{N}$, we have that
\begin{equation}
  \max_{-\rho<r<\rho}\left|u^{(1)}_\kappa-1\right|=o(\kappa^{-m}),
  \ \text{as}\ \kappa\to\infty.
\label{(16)}
\end{equation}
Next, we multiply $u^{(1)}_\kappa$ by a suitable constant of magnitude
$1$ and use (\ref{(12)}) and (\ref{(13)}) to introduce
$u^{(2)}_\kappa$ so that, in addition to (\ref{(14)}) and
(\ref{(16)}), it satisfies
\begin{equation}
{\rm Im}\int_{0}^{2\pi} u^{(2)}_\kappa(L,\varphi) d\varphi=
{\rm Im}\int_{0}^{2\pi} u^{(2)}_\kappa(-L,\varphi) d\varphi=0.
\label{(15)}
\end{equation}
Observe that 
\begin{multline*}
  {\left(\left|u^{(2)}_\kappa\right|^2-1\right)}^2=\left(\left({\rm
        Re}\left(u^{(2)}_\kappa\right)\right)^2+\left({\rm
        Im}\left(u^{(2)}_\kappa\right)\right)^2-1\right)^2 \\ \geq
  {\left({\rm Re}\left(u^{(2)}_\kappa\right)-1\right)}^2\left({\rm
      Re}\left(u^{(2)}_\kappa\right)+1\right)^2 -4\left(1-{\rm
      Re}\left(u^{(2)}_\kappa\right)\right)\left({\rm
      Im}\left(u^{(2)}_\kappa\right)\right)^2,
\end{multline*}
since $\left|u_\kappa\right|\leq 1$ by the maximum principle
\cite{BG2002} and, hence, ${\rm Re}\left(u^{(2)}_\kappa\right)\leq 1$.
% then, using (\ref{(14)}) and (\ref{(16)}), we have for sufficiently
% large $\kappa>0$ that there exists a $\kappa^\prime\geq\kappa$ such
% that $\tilde v_\kappa:=u_{\kappa^\prime}$ satisfies
% $$
% \tilde E_\kappa\left[\tilde v_\kappa\right]=\frac{1}{2}\int_{-L}^L
% dr\int_{0}^{2\pi}d\varphi |\nabla \tilde v_\kappa|^2
% $$
% $$+ \int_{-\rho}^\rho dr\int_{0}^{2\pi}d\varphi
% \left(\frac{\kappa^2}{2}({\rm Re}(\tilde v_\kappa)-1)^2
%   -\frac{\kappa^{-2}}{2}({\rm Im}(\tilde v_\kappa))^2 \right)\leq
% 2\pi.
% $$
  Then, using (\ref{(14)})--(\ref{(15)}) and ${\rm
    Re}\left(u^{(2)}_\kappa\right)\leq 1$ we have for any $m>0$ and
  any sufficiently large $\kappa>0$ that
\begin{eqnarray}
\label{ineq_suppl}  
L_\kappa\left[u^{(2)}_\kappa\right]&:=&\frac{1}{2}\int_{-L}^L dr\int_{0}^{2\pi}d\varphi \left|\nabla
  u^{(2)}_\kappa\right|^2 \\ &+&\int_{-\rho}^\rho dr\int_{0}^{2\pi}d\varphi
\left(\frac{\kappa^2}{4R^2}\left({\rm Re}\left(u^{(2)}_\kappa\right)-1\right)^2
    -{o\left(\kappa^{2-m}\right)}\left({\rm Im}\left(u^{(2)}_\kappa\right)\right)^2 \right)\leq 2\pi\,. \nonumber
\end{eqnarray}

Now let $m=5$ in (\ref{ineq_suppl}). Given a $\kappa>0$, we can choose
a sufficiently large $\kappa^\prime>\kappa$ such that
\begin{equation}
  \label{eq:1}
  \frac{\kappa^{\prime2}}{2R^2}>\kappa^2\ \mbox{and}\ {o\left(\kappa^{\prime-3}\right)}<\frac{\kappa^{-2}}{2}
\end{equation} 
and 
\begin{equation}
  \label{eq:2}
  L_\kappa\left[u^{(2)}_{\kappa^\prime}\right]\leq 2\pi\,.
\end{equation}
On the other hand, due to (\ref{eq:1}), we have that
\begin{equation}
  \label{eq:3}
  L_\kappa\left[w\right]\geq F_\kappa\left[w\right]:=\frac{1}{2}\int_{-L}^L
  dr\int_{0}^{2\pi}d\varphi |\nabla w|^2
  + \int_{-\rho}^\rho dr\int_{0}^{2\pi}d\varphi
  \left(\frac{\kappa^2}{2}({\rm Re}(w)-1)^2
    -\frac{\kappa^{-2}}{2}({\rm Im}(w))^2 \right)\,,
\end{equation}
for any function $w\in H^1\left((-L\,,L)\times(0\,,2\pi)\right)$. Note
that, unlike $E_\kappa[w]$, the functional $F_\kappa[w]$ is quadratic
in $w$ and, therefore, the Euler-Lagrange equation corresponding to
$F_\kappa$ is linear.

By substituting $v_\kappa:=u^{(2)}_{\kappa^\prime}$ in (\ref{eq:3})
and using (\ref{eq:2}), we obtain $F_\kappa\left[v_\kappa\right]\leq
2\pi$. Furthermore, $|v_\kappa|=1$ as $r=\pm L$, the function
$v_\kappa$ is $2\pi$-periodic in $\varphi$, and
$$
v_\kappa=a^\kappa_0+\sum_{n=1}^\infty(a^\kappa_n \cos
n\varphi +b^\kappa_n \sin n\varphi), \ \text{as}\ r=\pm L.
$$ 
In view of (\ref{(15)})
\begin{equation}
{\rm Im}(a^\kappa_0)=0
\label{im}
\end{equation} 
and
\begin{equation}
1=\frac{1}{2i}\sum_{n=1}^\infty n(b^\kappa_n\bar a^\kappa_n- a^\kappa_n\bar b^\kappa_n)=
\sum_{n=1}^\infty n({\rm Re}(a^\kappa_n){\rm Im}(b^\kappa_n)-{\rm Re}(b^\kappa_n){\rm Im}(a^\kappa_n)).
\label{(1)}
\end{equation}
by the degree formula $$ \mbox {\rm deg }(v,\Gamma)=\frac 1{2\pi
  i}\int_\Gamma \bar v\frac{\partial v}{\partial\tau},$$ valid when
$v\in C^1 (\Gamma ; S^1)$, $\Gamma$ is a $C^1$ simple closed curve in
$\D C$, and $\tau$ is a unit tangent vector to $\Gamma$ (cf.
\cite{bemi04a}).

For large $\kappa$ there is a unique minimizer $w_\kappa$ of
$F_\kappa\left[w\right]$ in the class of functions $2\pi$-periodic in
$\varphi$ and satisfying $w_\kappa=v_\kappa$ when $r=\pm L$.
Then
\begin{equation}
F_\kappa\left[w_\kappa\right]\leq F_\kappa\left[v_\kappa\right]\leq2\pi,
\label{(2)}
\end{equation}
where $w_\kappa$ is the solution of the problem
\begin{equation}
\label{linear}
\begin{cases}
  -\Delta {\rm Re}(w)+\kappa^2 V(r)({\rm Re}(w)-1)=0,\ -L<r<L,\\
  -\Delta {\rm Im}(w)-\kappa^{-2} V(r){\rm Im}(w)=0,\ -L<r<L,\\
  w(r,\varphi)=w(r,\varphi+2\pi),\\
  w=v_\kappa, \ r=\pm L.
\end{cases}
\end{equation}
Here $V(r)=1$ when $-\rho<r<\rho$ and $V(r)=0$ otherwise.  

\subsection{Energy estimate for the linear problem}
\label{estimate}
The problem \sy{linear} has the unique solution for large $\kappa$ in
the form
$$
w_\kappa(r,\varphi)=1+\left(a^\kappa_0-1\right) w^{(1)}_{\kappa,0}(r)+
\sum_{n=1}^\infty w^{(1)}_{\kappa,n}(r)({\rm
  Re}\left(a^\kappa_n\right)\cos{n\varphi}+ {\rm
  Re}\left(b^\kappa_n\right)\sin{n\varphi})
$$
$$
+i\sum_{n=1}^\infty w^{(2)}_{\kappa,n}(r)({\rm
  Im}\left(a^\kappa_n\right)\cos{n\varphi}+ {\rm
  Im}\left(b^\kappa_n\right)\sin{n\varphi})
$$
with real-valued $w^{(1)}_{\kappa,n}$ and $w^{(2)}_{\kappa,n}$ (here
it is important that $a^\kappa_0\in \D{R}$ by (\ref{im})).  The
functions $w^{(1)}_{\kappa,n}, w^{(2)}_{\kappa,n}$ can be found
explicitly so that
\begin{equation}
  F_\kappa\left[w_\kappa\right]=P^\kappa_0+\pi\sum_{n=1}^\infty n(P^\kappa_n(|{\rm
    Re}\left(a^\kappa_n\right)|^2+|{\rm Re}\left(b^\kappa_n\right)|^2)+
  Q^\kappa_n(|{\rm Im}\left(a^\kappa_n\right)|^2+|{\rm
    Im}\left(b^\kappa_n\right)|^2)).
\end{equation}
Here $P^\kappa_0\geq 0$ and the expressions for
$$
P^\kappa_n=
\frac{1-e^{-2n(L-\rho)}+\left(1+e^{-2n(L-\rho)}\right)\sqrt{1+\kappa^2
    n^{-2}}\tanh\left(\rho\sqrt{n^2+\kappa^2}\right)}
{1+e^{-2n(L-\rho)}+\left(1-e^{-2n(L-\rho)}\right)\sqrt{1+\kappa^2
    n^{-2}}\tanh\left(\rho\sqrt{n^2+\kappa^2}\right)},
$$
and
$$
Q^\kappa_n= \frac{1-e^{-2n(L-\rho)}+
  \left(1+e^{-2n(L-\rho)}\right)\sqrt{1-(\kappa n)^{-2}}
  \tanh\left(\rho\sqrt{n^2-\kappa^{-2}}\right)}
{1+e^{-2n(L-\rho)}+\left(1-e^{-2n(L-\rho)}\right)\sqrt{1-(\kappa
    n)^{-2}}\tanh\left(\rho\sqrt{n^2-\kappa^{-2}}\right)}.
$$
are derived in the Appendix. Then using $P^\kappa_0\geq 0$ and the elementary inequality $a^2+b^2>2ab$ we obtain
\begin{equation}
F_\kappa\left[w_\kappa\right]\geq 2\pi\sum_{n=1}^\infty n\sqrt{P^\kappa_n Q^\kappa_n} 
(|{\rm Re}\left(a^\kappa_n\right)||{\rm Im}\left(b^\kappa_n\right)|+|{\rm Re}\left(b^\kappa_n\right)||{\rm Im}\left(a^\kappa_n\right)|).
\label{(6)}
\end{equation}
Now we show that there exists a $\kappa_0>0$ such that
\begin{equation}
P^\kappa_nQ^\kappa_n>1,
\label{final}
\end{equation}
for all $\kappa\geq \kappa_0$ and all $n\geq 1$. Indeed, we can
rewrite $P^\kappa_n$ and $Q^\kappa_n$ as follows
$$
P^\kappa_n=
\frac{1+\beta^\kappa_n e^{-2n(L-\rho)}}{1-\beta^\kappa_n e^{-2n(L-\rho)}},\ \ 
Q^\kappa_n=
\frac{1-\alpha^\kappa_n e^{-2n(L-\rho)}}{1+\alpha^\kappa_n e^{-2n(L-\rho)}}
$$
where
$$
\alpha^\kappa_n= \frac{1-\sqrt{1-(\kappa
    n)^{-2}}\tanh\left(\rho\sqrt{n^2-\kappa^{-2}}\right)} {1+\sqrt{1-(\kappa
    n)^{-2}}\tanh\left(\rho\sqrt{n^2-\kappa^{-2}}\right)},
$$
and
$$
\beta^\kappa_n= \frac{\sqrt{1+\kappa^2
    n^{-2}}\tanh\left(\rho\sqrt{n^2+\kappa^2}\right)-1}
{\sqrt{1+\kappa^2
    n^{-2}}\tanh\left(\rho\sqrt{n^2+\kappa^2}\right)+1}.
$$
Note that (\ref{final}) is equivalent to the inequality
$\alpha^\kappa_n<\beta^\kappa_n$.  This inequality clearly holds for
any fixed $n\geq 0$ when $\kappa$ is sufficiently large, since
$$
\alpha^\kappa_n\to e^{-2n\rho},\ \ \beta^\kappa_n\to 1, \ \ \text{as} \ \ \kappa\to\infty.
$$
On the other hand, for all $\kappa\geq 1$, multiplying and
dividing $\alpha^\kappa_n$ and $\beta^\kappa_n$ by their respective
denominators and letting $n\to \infty$, we have
$$
\alpha^\kappa_n\leq e^{-n\rho}+\frac{1}{(n\kappa)^{2}},\ \
\beta^\kappa_n\geq\frac{\gamma}{n^2},
$$
where $\gamma>0$ is independent of $n$ and $\kappa$. Thus
$\alpha^\kappa_n<\beta^\kappa_n$ and, hence, (\ref{final}) are
satisfied once $\kappa_0$ is chosen to be sufficiently large.

By (\ref{(6)}) and (\ref{final}) we get
$$
F_\kappa\left[w_\kappa\right]\geq
2\pi\sum_{n=1}^\infty n(|{\rm Re}\left(a^\kappa_n\right)||{\rm
  Im}\left(b^\kappa_n\right)|+|{\rm Re}\left(b^\kappa_n\right)||{\rm
  Im}\left(a^\kappa_n\right)|),
$$
and, according to (\ref{final}), this inequality is strict unless
r.h.s.=0. By (\ref{(1)})
$$
\sum_{n=1}^\infty n\left(\left|{\rm
      Re}\left(a^\kappa_n\right)\right|\left|{\rm
      Im}\left(b^\kappa_n\right)\right|+\left|{\rm
      Re}\left(b^\kappa_n\right)\right|\left|{\rm
      Im}\left(a^\kappa_n\right)\right|\right)\geq 1,
$$
so that $F_\kappa\left[w_\kappa\right]>2\pi$. This
contradicts (\ref{(2)}).

\bibliographystyle{plain}
\bibliography{nonexistence}

\section{Appendix. Computation of $P^\kappa_n$ and $Q^\kappa_n$.}

Suppose that $w_\kappa$ is the solution of \sy{linear}. Multiply the
first equation in \sy{linear} by ${\mathrm{Re}}(w_\kappa)-1$ and the
second equation in \sy{linear} by ${\mathrm{Im}}(w_\kappa)$. Adding
the resulting equations together, integrating by parts, and using the
symmetry of $w_\kappa$ we obtain
\begin{multline}
  F_\kappa\left[ w_\kappa\right]=\frac{1}{2}\int^{2\pi}_{0}
  \left(\left({\rm
        Re}\left(w_\kappa\right)\right)(L,\varphi)-1\right)\frac{d({\rm Re}\left(w_\kappa\right))}{dr}(L,\varphi)\, d\varphi \\
  -\frac{1}{2}\int^{2\pi}_{0} \left(\left({\rm Re}\left(
        w_\kappa\right)\right)(-L,\varphi)-1\right)\frac{d({\rm
      Re}\left(
      w_\kappa\right))}{dr}(-L,\varphi)\, d\varphi \\
  +\frac{1}{2}\int^{2\pi}_{0} ({\rm Im}\left(
    w_\kappa\right))(L,\varphi)\frac{d({\rm
      Im}\left(w_\kappa\right))}{dr}(L,\varphi)\, d\varphi \\
  -\frac{1}{2}\int^{2\pi}_{0} ({\rm Im}\left(
    w_\kappa\right))(-L,\varphi)\frac{d({\rm Im}\left(
      w_\kappa\right))}{dr}(-L,\varphi)\, d\varphi, \\
  =\int^{2\pi}_{0} \left(\left({\rm
        Re}\left(w_\kappa\right)\right)(L,\varphi)-1\right)\frac{d({\rm
      Re}\left(w_\kappa\right))}{dr}(L,\varphi)+({\rm Im}\left(
    w_\kappa\right))(L,\varphi)\frac{d({\rm
      Im}\left(w_\kappa\right))}{dr}(L,\varphi)\, d\varphi
\label{enparts}
\end{multline}
Further, substituting the expansions
\begin{equation}
  \label{eq:a1}
  {\rm Re}\left(w_\kappa\right)=a^\kappa_0
w^{(1)}_{\kappa,0}(r)+ \sum_{n=1}^\infty w^{(1)}_{\kappa,n}(r)({\rm
  Re}\left(a^\kappa_n\right)\cos{n\varphi}+ {\rm
  Re}\left(b^\kappa_n\right)\sin{n\varphi}),  
\end{equation}
\begin{equation}
  \label{eq:a2}
  {\rm Im}\left(w_\kappa\right)=\sum_{n=1}^\infty
w^{(2)}_{\kappa,n}(r)({\rm Im}\left(a^\kappa_n\right)\cos{n\varphi}+
{\rm Im}\left(b^\kappa_n\right)\sin{n\varphi}),  
\end{equation}
into \sy{enparts} and integrating, the expression for
$F_\kappa[w_\kappa]$ can be written as
\begin{multline}
  F_\kappa\left[ w_\kappa\right]=a_0^\kappa \frac{d}{dr}w^{(1)}_{\kappa,0}(L)\left(a_0^\kappa w^{(1)}_{\kappa,0}(L)-1\right) \\
  +\pi\sum_{n=1}^\infty
  \left(w^{(1)}_{\kappa,n}(L)\frac{d}{dr}w^{(1)}_{\kappa,n}(L)\left(|{\rm
        Re}\left(a^\kappa_n\right)|^2+|{\rm Re}\left(b^\kappa_n\right)|^2\right)\right. \\
  \left.+w^{(2)}_{\kappa,n}(L)\frac{d}{dr}w^{(2)}_{\kappa,n}(L)\left(|{\rm
        Im}\left(a^\kappa_n\right)|^2+|{\rm
        Im}\left(b^\kappa_n\right)|^2\right)\right).
\label{eq:a6}
\end{multline}
We set
\begin{equation}
  \label{eq:a7}
  P^\kappa_0:=a_0^\kappa
  \frac{d}{dr}w^{(1)}_{\kappa,0}(L)\left(a_0^\kappa
    w^{(1)}_{\kappa,0}(L)-1\right)\,,
\end{equation}
and
\begin{equation}
  \label{eq:a8}
  P^\kappa_n:=w^{(1)}_{\kappa,n}(L)\frac{d}{dr}w^{(1)}_{\kappa,n}(L)\,,\ \ Q_n^\kappa:=w^{(2)}_{\kappa,n}(L)\frac{d}{dr}w^{(2)}_{\kappa,n}(L)\,,
\end{equation}
for every $n\geq 1$.

If we assume that ${\rm Re}(w_\kappa)=const$ and ${\rm
  Im}(w_\kappa)=0$ on $\left\{-L\,,L\right\}\times[0\,,2\pi]$, then
$a^\kappa_n=0$ and $b^\kappa_n=0$ for all $n\geq 1$. The remaining
term in \sy{eq:a6} must be nonnegative because $F_\kappa[w_\kappa]\geq
0$ when ${\rm Im}(w_\kappa)\equiv0$.  We conclude that $P^\kappa_0\geq
0$.

Using the standard separation of variables argument, we have that the
functions $w^{(1)}_{\kappa,n}$ and $w^{(2)}_{\kappa,n}$ satisfy
\begin{equation}
  \label{eq:a3}
\begin{cases}
  \displaystyle -\frac{d^2}{dr^2}w^{(1)}_{\kappa,n}(r)+\left(n^2+\kappa^2 V(r)\right)w^{(1)}_{\kappa,n}(r)=0, & -L<r<L,\\ \\
  w^{(1)}_{\kappa,n}(\pm L)=1, &
\end{cases}  
\end{equation}
and
\begin{equation}
  \label{eq:a4}
  \ \begin{cases}
    \displaystyle -\frac{d^2}{dr^2}w^{(2)}_{\kappa,n}(r)+\left(n^2-\kappa^{-2} V(r)\right)w^{(2)}_{\kappa,n}(r)=0, & -L<r<L,\\ \\
    w^{(2)}_{\kappa,n}(\pm L)=1, &
\end{cases}  
\end{equation}
for $n\geq1$. 

Solving \sy{eq:a3} and \sy{eq:a4}, we obtain
\[
\gamma^{(1)}_{\kappa,n}\,w^{(1)}_{\kappa,n}(r) =
\begin{cases}
  \begin{aligned}\cosh&\,{\left(n(r-\rho)\right)}\cosh{\left(\rho\sqrt{n^2+\kappa^2}\right)}
    \\+&\frac{\sqrt{n^2+\kappa^2}}{n}\sinh{\left(n(r-\rho)\right)}\sinh{\left(\rho\sqrt{n^2+\kappa^2}\right)},\end{aligned}& \text{if }r\in(\rho,L),\\ \\
  \,\cosh{\left(r\sqrt{n^2+\kappa^2}\right)},& \text{if }r\in(-\rho,\rho),\\ \\
  \begin{aligned}\cosh&\,{\left(n(r+\rho)\right)}\cosh{\left(\rho\sqrt{n^2+\kappa^2}\right)} 
    \\-&\frac{\sqrt{n^2+\kappa^2}}{n}\sinh{\left(n(r+\rho)\right)}\sinh{\left(\rho\sqrt{n^2+\kappa^2}\right)},\end{aligned}&
  \text{if }r\in(-L,-\rho),
\end{cases}
\]
\[
  \gamma^{(2)}_{\kappa,n}\,w^{(2)}_{\kappa,n}(r) =
\begin{cases}
  \begin{aligned}\cosh&\,{\left(n(r-\rho)\right)}\cosh{\left(\rho\sqrt{n^2-\kappa^{-2}}\right)}
    \\+&\frac{\sqrt{n^2-\kappa^{-2}}}{n}\sinh{\left(n(r-\rho)\right)}\sinh{\left(\rho\sqrt{n^2-\kappa^{-2}}\right)},\end{aligned}& \text{if }r\in(\rho,L),\\ \\
  \,\cosh{\sqrt{n^2-\kappa^{-2}}r},& \text{if }r\in(-\rho,\rho),\\ \\
  \begin{aligned}\cosh&\,{\left(n(r+\rho)\right)}\cosh{\left(\rho\sqrt{n^2-\kappa^{-2}}\right)} 
    \\-&\frac{\sqrt{n^2-\kappa^{-2}}}{n}\sinh{\left(n(r+\rho)\right)}\sinh{\left(\rho\sqrt{n^2-\kappa^{-2}}\right)},\end{aligned}&
  \text{if }r\in(-L,-\rho),
\end{cases}
\]
where
\[
\gamma^{(1)}_{\kappa,n}=\cosh{\left(n(L-\rho)\right)} \
\cosh{\left(\rho\sqrt{n^2+\kappa^2}\right)}
+\frac{\sqrt{n^2+\kappa^2}}{n} \sinh{\left(n(L-\rho)\right)} \ \sinh
\sqrt{n^2+\kappa^2}\rho,
\]
\[
\gamma^{(2)}_{\kappa,n}= \cosh{\left(n(L-\rho)\right)}\,\cosh
{\left(\rho\sqrt{n^2-\kappa^{-2}}\right)} +\frac{\sqrt{n^2-\kappa^{-2}}}{n} \sinh{\left(
n(L-\rho)\right)}\,\sinh{\left(\rho\sqrt{n^2-\kappa^{-2}}\right)}.
\]
Substituting the expressions for $w^{(1)}_{\kappa,n}$ and
$w^{(2)}_{\kappa,n}$ into \sy{eq:a8}, we have
$$
P^\kappa_n=
\frac{1-e^{-2n(L-\rho)}+\left(1+e^{-2n(L-\rho)}\right)\sqrt{1+\kappa^2
    n^{-2}}\tanh\left(\rho\sqrt{n^2+\kappa^2}\right)}
{1+e^{-2n(L-\rho)}+\left(1-e^{-2n(L-\rho)}\right)\sqrt{1+\kappa^2
    n^{-2}}\tanh\left(\rho\sqrt{n^2+\kappa^2}\right)},
$$
and
$$
Q^\kappa_n= \frac{1-e^{-2n(L-\rho)}+
  \left(1+e^{-2n(L-\rho)}\right)\sqrt{1-(\kappa n)^{-2}}
  \tanh\left(\rho\sqrt{n^2-\kappa^{-2}}\right)}
{1+e^{-2n(L-\rho)}+\left(1-e^{-2n(L-\rho)}\right)\sqrt{1-(\kappa
    n)^{-2}}\tanh\left(\rho\sqrt{n^2-\kappa^{-2}}\right)}.
$$
\end{document}